\documentclass[11pt]{article}

\usepackage[utf8]{inputenc} 
\usepackage[T2A,OT1]{fontenc} 
\usepackage[english,russian]{babel}

\newcounter{isappendix}
\def\appendix%
{\section*{\hfill\sl П\,Р\,И\,Л\,О\,Ж\,Е\,Н\,И\,Е}%
	\refstepcounter{isappendix}}

\newtheorem{lemma}{\hspace{\parindent}\sl{Л\,е\,м\,м\,а\,}}

\newcommand*{\No}{\textnumero}

\usepackage{graphicx}
\usepackage{color}
\usepackage{cmap}
\usepackage{amsmath,amsfonts,amssymb,
amsbsy,amstext,amscd,amsxtra,multicol}
\usepackage{units}
\usepackage{fancyhdr}
\usepackage{forloop}
\usepackage{indentfirst}
\usepackage{verbatim}
\usepackage[colorlinks, urlcolor=blue]{hyperref}
\usepackage{url}
\usepackage{algorithm}
\usepackage{algpseudocode}
\usepackage{graphicx}
\usepackage{pifont}
\usepackage[inline]{asymptote}

\makeatletter
\renewcommand{\fnum@algorithm}{\fname@algorithm}
\makeatother

\def\Prog{\text{Prog}}

\newcommand{\E}{\mathbb{E}}
\newcommand{\RR}{\mathbb{R}}

\newcommand{\argmin}{\mathop{\arg\!\min}}

\newcommand{\circledOne}{\text{\ding{172}}}
\newcommand{\circledTwo}{\text{\ding{173}}}
\newcommand{\circledThree}{\text{\ding{174}}}
\newcommand{\circledFour}{\text{\ding{175}}}

\newcounter{fakecnt}

\newtheorem{LemAp}{\hspace{\parindent}{\sl Л~е~м~м~а~}}[fakecnt]

\newtheorem{CoroNoNum}{\hspace{\parindent}{\sl У~т~в~е~р~ж~д~е~н~и~е~}}

\newtheorem{TheoremNoNum}{\hspace{\parindent}{\sl Т~е~о~р~е~м~а~}}

\title{%
	Ускоренный спуск по случайному направлению 
	с неевклидовой прокс-структурой 
\footnote{Работа А.В.~Гасникова по разделам 3 и 4 поддержана Российским фондом фундаментальных исследований (проект \No~18-31-20005 мол\_а\_вед). Работа Э.А.~Горбунова и Е.А.~Воронцовой поддержана грантом Президента РФ~МД-1320.2018.1. Работа А.В.~Гасникова и  Е.А.~Воронцовой поддержана Российским фондом фундаментальных исследований (проект \No~18-29-03071 мк).}%
}%

\author{Е.А.~ВОРОНЦОВА, канд.~ф.-м.~наук (\textrm{vorontsovaea@gmail.com}) \\
(Дальневосточный федеральный университет, Владивосток),
	\\
	А.В.~ГАСНИКОВ, д-р~ф.-м.~наук (\textrm{gasnikov@yandex.ru}), \\
	Э.А.~ГОРБУНОВ (\textrm{ed-gorbunov@yandex.ru})
	\\
	(Московский физико-технический институт) }

\date{}

\begin{document}
	
	\maketitle
\small{$\quad$Рассматриваются задачи гладкой выпуклой оптимизации,.
для численного решения которых полный градиент недоступен.
В 2011 г. Ю.Е.~Нестеровым были предложены ускоренные безградиентные методы решения таких задач. 
Поскольку рассматривались только задачи безусловной оптимизации, то использовалась евклидова прокс-структура. Однако если заранее знать, например, что решение задачи разреженно, а точнее, что расстояние от точки старта до решения в 1-норме и в 2-норме близки, то более выгодно выбирать не евклидову прокс-структуру, связанную с 2-нормой, а прокс-структуру, связанную с 1-нормой. 
Полное обоснование этого утверждения проводится в статье.
Предлагается
ускоренный метод спуска по случайному направлению с неевклидовой прокс-структурой для решения задачи безусловной оптимизации (в дальнейшем подход
предполагается расширить на ускоренный безградиентный метод).
Получены оценки скорости сходимости метода.
Показаны сложности переноса 
описанного подхода на задачи условной оптимизации.
}

	\textit{Ключевые слова:} ускоренные методы первого порядка, выпуклая оптимизация, метод линейного каплинга, концентрация равномерной меры на единичной евклидовой сфере, неевклидова прокс-структура.

	\section{Введение}
	В \cite{nest} были предложены ускоренные оракульные \footnote{Здесь и далее под 
	оракулом понимается подпрограмма расчета
	значений целевой функции и/или градиента (его части),
	а 	оптимальность метода на классе задач понимается в смысле Бахвалова--Немировского~\cite{Nemir} как число обращений (по ходу работы метода) к оракулу для достижения заданной точности (по функции). } $\,$ методы нулевого порядка (безградиентные методы) решения задач гладкой выпуклой безусловной оптимизации.
	
	В рассуждениях \cite{nest} существенным образом использовалось то, что была выбрана евклидова прокс-структура (выпуклая гладкая функция, порождающая расстояние, и 1-сильно выпуклая относительно какой-то нормы (строгое определение см. в разделе~4)). Такой выбор прокс-структуры для задач безусловной оптимизации является вполне естественным (см., например, \cite{stochastic-gradients-inexact-oracle}). Однако в ряде задач имеется дополнительная информация, которая, например, позволяет рассчитывать на разреженность решения (в решении большая часть компонент нулевые). В таких случаях использование других прокс-структур бывает более выгодным. Для негладких задач стохастической условной оптимизации с оракулом нулевого порядка недавно было показано (см. \cite{GLUF, GKLUF}), что в определенных ситуациях ускорение метода за счет перехода от евклидовой прокс-структуры, связанной с 2-нормой, к прокс-структурам, связанным с 1-нормой, может давать ускорение методу, по порядку равное размерности пространства, в котором происходит оптимизация. К сожалению, техника, использованная в \cite{GLUF, GKLUF} существенным образом использовала неускоренную природу оптимальных методов для негладких задач. Другими словами, из \cite{GLUF, GKLUF} непонятно, как получать аналогичные оценки для гладких задач. В настоящей статье на базе специального варианта быстрого (ускоренного) градиентного метода \cite{Allen-Zhu-Orecchia-linear-coupling} строится ускоренный метод спуска по случайному направлению.
	Особенностью метода из \cite{Allen-Zhu-Orecchia-linear-coupling} является представление быстрого градиентного метода как специальной выпуклой комбинации градиентного спуска и зеркального спуска. В \cite{Allen-Zhu-Orecchia-linear-coupling}, как и во всех известных авторам вариантах быстрого градиентного метода с двумя и более $\,$ ``проекциями'', обе проекции осуществлялись в одной норме/прокс-структуре. Главной идеей настоящей статьи является использование разных норм/прокс-структур в этих проекциях, а именно: в градиентном шаге всегда используется обычная евклидова проекция, а вот в зеркальном шаге выбор прокс-структуры обусловлен априорной информацией о свойствах решения.
	
	В настоящей статье на базе описанной конструкции для детерминированных задач безусловной гладкой выпуклой оптимизации строится ускоренный метод спуска по случайному направлению \footnote{Подробно о разнице 
	в подходах в случае детерминированной постановки,
	но с введением рандомизации и в случае
	задач стохастической оптимизации
	см. в~\cite{Gasn_Universal_2018}.}  (раздел~4).
	
	В классе детерминированных спусков по направлению (к ним можно отнести и циклический координатный спуск)
	для получения лучших оценок необходимо
	вводить рандомизацию (доказательство см. в.~\cite{Nest_coord_des}), поэтому
	рассматриваются сразу спуски по случайному направлению.

	Показано, какие возникают сложности при попытке перенесения описанного подхода на задачи оптимизации на множествах простой структуры (раздел~4). 

	\section{Постановка задачи}
	Рассматривается задача гладкой выпуклой оптимизации
	\begin{equation} \label{optimization_problem}
		f\left( x \right)\to \mathop {\min }\limits_{x\in Q},
	\end{equation}
	где функция $f\left( x \right)$, заданная на выпуклом замкнутом множестве $Q \subseteq \RR^n$, имеет липшицев градиент с 
	константой $L_2$ (т.е. $f(x)$~--- $L_2$-гладкая функция)
	\[
	\left\| {\nabla f\left( y \right)-\nabla f\left( x \right)} \right\|_2 \le 
	L_2 \left\| {y-x} \right\|_2
	\]
	и является $\mu $-сильно выпуклой в $p$-норме ($1\le p\le 2)$ функцией (далее будем использовать и обозначение $\mu_p$), при этом в точке минимума $x_*$ выполнено равенство $\nabla f\left( {x_\ast } \right)=0$, а итерационный процесс стартует с точки~$x_0$. 
	
	В данной статье для решения задачи~\eqref{optimization_problem} вместо обычного градиента используется его стохастическая
	аппроксимация, построенная на базе производной по случайно выбранному направлению 
    \cite{triangles-method}
	\[
	g\left( {x, \, e} \right)\overset{\text{def}}{=}n\left\langle {\nabla f\left( x \right), \, e} 
	\right\rangle e,
	\]
	где $e$~--- случайный вектор, равномерно распределенный на $S_2^n \left( 1 
	\right)$~--- единичной сфере в 2-норме в пространстве $\RR^n$ ($e \sim RS_2^n \left( 1 \right)$; под этой записью будем понимать, что случайный вектор $e$ имеет равномерное распределение на $n$-мерной единичной евклидовой сфере), а угловые скобки~$\langle\cdot, \, \cdot\rangle$ обозначают скалярное произведение \footnote{Отметим, что $\E_e[g(x,e)] = \nabla f(x),$ что следует из факта: $\E_e[|e_i|^2] = \frac{1}{n}$, где $e_i$~--- $i$-я компонента вектора $e$.}.

	Имеет место следующая лемма (доказательство см. в~\cite{GorVorGas_lemma1}), являющаяся следствием явления 
	концентрации равномерной меры на сфере вокруг экватора (см. также~\cite{BD}; северный полюс задается градиентом $\nabla f\left( x \right)$).
	\begin{lemma}
		Пусть $e \sim RS_2^n(1),\, n \geqslant8,\, s\in\RR^n$, тогда
		\begin{equation}\label{lemm1:expect_q_norm}
        \E[||e||_q^2] \leqslant \min\{q-1,\,16\ln n - 8\}n^{\frac{2}{q}-1},\quad 2\leqslant q \leqslant \infty,
    \end{equation}
    
    \begin{equation}\label{lemm1:expect_inner_product}
        \E[\langle s,\, e\rangle^2||e||_q^2] \leqslant \sqrt{3}||s||_2^2\min\{2q-1,32\ln n -8\}n^{\frac{2}{q}-2},\quad 2\leqslant q \leqslant \infty,
    \end{equation}
    где под знаком $||\cdot||_q$ понимается векторная~$q$-норма $($норма Гельдера с показателем~$q)$.
	\end{lemma}
	
	Из \eqref{lemm1:expect_inner_product} (см. также~\cite{GLUF}) вытекает следующий факт.
	
	\begin{CoroNoNum}
		Пусть $e \sim RS_2^n \left( 1 \right)$ и $g\left( {x, \, e} \right)=n\left\langle {\nabla f\left( x \right), \, e} \right\rangle e$, тогда
		\[
		\E_e \left[ {\left\| {g\left( {x, \, e} \right)} \right\|_q^2 } \right]\leqslant\sqrt{3} {\min\{2q-1,\,32\ln n -     8\} n^{\frac{2}{q}}\left\| {\nabla f\left( x \right)} \right\|_2^2 
		} .
				\]
	\end{CoroNoNum}

	Используя данное утверждение, можно в сильно выпуклом случае при условии $\nabla 
	f\left( {x_\ast } \right)=0$ получить оценку необходимого числа обращений к оракулу за производной по направлению для достижения по функции точности $\varepsilon$ в среднем~\cite{about-non-triviality}:
	\[
	N\left( \varepsilon \right)= {{\rm O}}\left( {n^{\frac{2}{q}} \ln n  \frac{L_2 
		}{\mu }\ln \left( {\frac{\Delta f^0}{\varepsilon }} \right)} \right),
	\]
	где $\Delta f^0 = f(x_0) - f(x_{*})$.
	
	Имеется гипотеза (см., например, \cite{stochastic-gradients-inexact-oracle}), 
	что, используя специальные ускоренные методы (типа Катюши (Katusha) из 
	\cite{Allen-Zhu-Katusha}), можно получить оценку
	\[
	N\left( \varepsilon \right)={\rm O}\left( {n^{\frac{1}{q}+\frac{1}{2}}\ln n\sqrt 
		{\frac{L_2 }{\mu }} \ln \left( {\frac{\Delta f^0}{\varepsilon }} \right)} 
	\right).
	\]
	Насколько известно авторам, эта гипотеза на данный момент не доказана и не опровергнута.
	
	В разделе 4 данной статьи для случая $Q = \RR^n$ доказывается оценка
	\begin{equation}\label{main_estimation}
		\E[f\left( {x_N } \right)] - f(x_{*}) \le Cn^{\frac{2}{q}+1}\ln n\frac{L_2 R^2}{N^2}.
	\end{equation}
	Эта оценка получается из приведенной выше оценки с помощью регуляризации $\mu = \varepsilon/R^2$ \cite{Gasn_Universal_2018}, где (с точностью до корня из логарифмического по $n$ множителя) $R$~--- расстояние в $p$-норме от точки старта до решения.

	\section{Задача А.С.~Немировского}
	Рассмотрим задачу \eqref{optimization_problem} минимизации гладкого выпуклого функционала 
	$f\left( x \right)$ с константой Липшица градиента $L_2 $ в 2-норме  
	на множестве $Q=B_1^n \left( R \right)$ (шар в пространстве ${\RR}^n$ 
	радиуса $R$ в 1-норме). Тогда 
	 на рассматриваемом классе функции для любого итерационного метода, на каждой итерации которого только один раз можно обратиться к оракулу за градиентом функции, можно так подобрать функцию $f(x)$ из этого класса, что имеет место оценка скорости сходимости \cite{guzman-nemirovski} в виде
	
	\begin{equation}\label{lower_bound}
		f\left( {x_N} \right)-f(x_*) \ge \frac{\tilde {C}_1 L_2 R^2}{N^3},
	\end{equation}
	где $\tilde{C}_1$~--- некоторая числовая константа \footnote{Как и $\tilde{C}_2$ далее. Здесь и далее все числовые константы не зависят от $N$ и $n$.}, $x^*$~--- ближайшая к $x_0$ точка минимума функции $f(x)$. Поскольку
	$f\left( {x_N} \right)-f(x_*)$ должно быть меньше или
	равно $\varepsilon$, из~\eqref{lower_bound} можно получить оценку на $N(\varepsilon)$.
	
	С другой стороны, если использовать обычный быстрый градиентный метод с 
	KL-прокс-структурой 
	для этой же задачи, то 
	\cite{universal-composite}:
	\[
	f\left( {x_N} \right)-f(x_*) \le \frac{\tilde {C}_2 L_1 R^2\ln n}{N^2},
	\]
	где константа Липшица градиента $L_1$ в 1-норме  удовлетворяет условию: ${L_2 } \mathord{\left/ {\vphantom 
			{{L_2 } n}} \right. \kern-\nulldelimiterspace} n\le L_1 \le L_2 \,$ \footnote{Здесь под $L_1$ понимается такое положительное число, что $||\nabla f(x) - \nabla f(y)||_\infty \leqslant L_1||x-y||_1$.}. 
	Отсюда нельзя сделать вывод, что нижняя оценка достигается. 
	Достигается ли эта нижняя оценка и если достигается, то на каком методе? Насколько авторам известно, это пока открытый вопрос, поставленный А.С.~Немировским в 2015 г. (см. также~\cite{guzman-nemirovski}). Однако если оценивать не число итераций, а общее 
	число арифметических операций и если ограничиться рассмотрением класса 
	функций, для которых ``стоимость'' расчета производной по направлению примерно в $n$ 
	раз меньше ``стоимости'' расчета полного градиента \footnote{Из-за быстрого 
	автоматического дифференцирования \cite{AD} это предположение 
	довольно обременительное; но если функция задана моделью черного 
	ящика, выдающего только значение функции, а градиент восстанавливается при 
	$n+1$ таком обращении, то сделанное предположение кажется вполне 
	естественным.}, то при $N \le n$ (точнее даже при $N\simeq n)$ выписанная оценка \eqref{main_estimation} (в варианте для общего числа арифметических операций, 
	необходимых для достижения заданной точности в среднем) с точностью до 
	логарифмического множителя будет соответствовать нижней оценке~\eqref{lower_bound}, если последнюю
	понимать как
	$$
	f\left( {x_N} \right)-f(x_*) \ge \frac{\tilde {C}_1 L_2 R^2}{nN^2},
	$$
	т.е.
	$$N\left( \varepsilon \right)={\rm O}\left( {\sqrt {\frac{L_2 
				R^2}{\varepsilon n}} } \right).$$

	Действительно ($q=\infty$ ), общее число арифметических операций равно
	\[
	{\rm O}\left( n \right)\cdot \underbrace { {{\rm O}}\left( 
		{n^{\frac{1}{2}}\ln n  \sqrt {\frac{L_2 R^2}{\varepsilon }} } 
		\right)}_{\text{число\;итераций}} \, \approx \, {\rm O}\left( {n^2} \right)\cdot 
	{{\rm O}\left( {\sqrt {\frac{L_2 R^2}{\varepsilon n}} } 
		\right)}.
	\]
	
	\section{Обоснование формулы \eqref{main_estimation} в случае $Q = \RR^n$}
	Введем дивергенцию Брэгмана~\cite{Bregman} $V_z \left( y \right)$, связанную с $p$-нормой ($1\leqslant p \leqslant 2$)
	$$
	V_z(y) \overset{\text{def}}{=} d(y) - d(z) - \langle\nabla d(z), \, y-z\rangle,
	$$
	где функция $d(x)$ является непрерывно дифференцируемой сильно и выпуклой с константой сильной выпуклости, равной единице. Например, для $p=1$ функцию $d(x)$ можно выбрать так:
	$$
	d(x) = \frac{1}{2(a-1)}||x||_a^2,
	$$
	где $a = \frac{2\log n}{2\log n - 1}$. Функцию $d(x)$ будем называть \textit{прокс-функцией} (или \textit{прокс-структурой}), связанной с $p$-нормой. Кроме того, пусть $q$~--- такое число, что $\frac{1}{p}+\frac{1}{q}=1$.  Далее будем следовать обозначениям из \cite{about-non-triviality}. Пусть случайный вектор $e$ равномерно распределен на поверхности евклидовой сферы единичного 
	радиуса ($e\sim RS_2^n \left( 1 \right))$. Положим, что
	\[
	\mbox{Grad}_e \left( x \right)\overset{\text{def}}{=}x-\frac{1}{L }\left\langle {\nabla f\left( 
		x \right), \, e} \right\rangle e,
	\]
	\[
	\mbox{Mirr}_e \left( {x,z,\alpha } \right)\overset{\text{def}}{=} \argmin\limits_{y\in\RR^n} \left\{ {\alpha \left\langle {n\left\langle {\nabla f\left( x 
				\right), \, e} \right\rangle e, \, y-z} \right\rangle +V_z \left( y \right)} 
	\right\},
	\]
	где $L$~--- константа Липшица градиента функции $f(x)$ в $2$-норме (индекс 2 не пишем, так как везде далее интересуемся константой Липшица в 2-норме).
	
	Опишем ускоренный неевклидов спуск (английское название метода~--- Accelerated by Coupling Directional Search, ACDS), построенный на базе специальной комбинации спусков по 
	направлению в форме градиентного спуска (Grad) и метода зеркального спуска (Mirr).
	\floatname{algorithm}{Алгоритм 1.}
	\begin{algorithm}
		\caption{Ускоренный неевклидов спуск (ACDS)}\label{ACRCD}
		\begin{algorithmic}[1]
			
			\Require $f$~--- выпуклая дифференцируемая функция на $\RR^n$ с липшицевым градиентом с константой $L$ по отношению к $2$-норме; $x_0$~--- некоторая стартовая точка; $N$~--- количество итераций.
			\Ensure точка $y_N$, для которой выполняется $\E_{e_1,e_2,\ldots,e_N}[f(y_N)] - f(x^*) \leqslant \frac{4\Theta LC_{n,q}}{N^2}$.
			\State $y_0 \leftarrow x_0, \, z_0 \leftarrow x_0$
			\For{$k=0,\, \dots, \, N-1$}
			\State $\alpha_{k+1} \leftarrow \frac{k+2}{2LC_{n,q}}, \, \tau_{k} \leftarrow \frac{1}{\alpha_{k+1}LC_{n,q}} = \frac{2}{k+2}$
			\State Генерируется $e_{k+1} \sim RS_2^n \left( 1 \right)$ независимо от предыдущих итераций
			\State $x_{k+1} \leftarrow \tau_kz_k + (1-\tau_k)y_k $
			\State $y_{k+1} \leftarrow \text{Grad}_{e_{k+1}}(x_{k+1})$
			\State $z_{k+1} \leftarrow \text{Mirr}_{e_{k+1}}(x_{k+1}, \, z_{k}, \, \alpha_{k+1})$
			\EndFor
			\State\Return $y_N$
		\end{algorithmic}
	\end{algorithm}

	\begin{TheoremNoNum}
		Пусть $f(x)$~--- выпуклая дифференцируемая функция на $Q=\RR^n$ с константой Липшица для градиента, равной $L$ в 2-норме, $d(x)$~--- 1-сильно выпуклая в $p$-норме функция на $Q$, $N$ -- число итераций метода, $x^*$~--- точка минимума функции $f(x)$. Тогда ускоренный неевклидов спуск $($ACDS$)$ на выходе даст точку $y_N$, удовлетворяющую неравенству
		$$
		\E_{e_1,e_2,\ldots,e_N}[f(y_N)] - f(x^*) \leqslant \frac{4\Theta LC_{n,q}}{N^2},
		$$
		где $\Theta \overset{\text{def}}{=} V_{x_0}(x^*)$, $C_{n,q} \overset{\text{def}}{=} \sqrt{3}\min\{2q-1,32\ln n -8\}n^{\frac{2}{q}+1}$,  $\frac{1}{q}+\frac{1}{p}=1$.
	\end{TheoremNoNum}
	
	Сформулируем две ключевые леммы, которые понадобятся для доказательства теоремы 
	(доказательства приведены в Приложении).
	
	\begin{lemma}
		Если $\tau_{k} = \frac{1}{\alpha_{k+1}LC_{n,q}}$, то для всех $u\in Q = \RR^n$ верны неравенства
		\begin{equation}\label{grad_estimation}
					\begin{array}{l}
				\alpha_{k+1}\langle \nabla f(x_{k+1}), z_{k}-u\rangle \leqslant\\ 
				\leqslant \alpha_{k+1}^2\cdot\frac{C_{n,q}}{2n}||\nabla f(x_{k+1})||_2^2 + V_{z_k}(u) \,-\\
				- \E_{e_{k+1}}[V_{z_{k+1}}(u)\mid e_1, \, e_2, \, \ldots, \, e_k] \leqslant\\ \leqslant \alpha_{k+1}^2LC_{n,q}\cdot(f(x_{k+1}) - \E_{e_{k+1}}[f(y_{k+1})\mid e_1, \, e_2, \, \ldots, \, e_k]) \, +\\
				+ \, V_{z_k}(u) - \E_{e_{k+1}}[V_{z_{k+1}}(u)\mid e_1, \, e_2, \, \ldots, \, e_k].
			\end{array}
		\end{equation}
	\end{lemma}
	
	\begin{lemma}
		Для всех $u \in Q=\RR^n$ выполнено неравенство
		\begin{equation}\label{coupling}
			\begin{array}{l}
				\alpha_{k+1}^2 L C_n\E_{e_{k+1}}[f(y_{k+1})\mid e_1,\ldots,e_{k}] - (\alpha_{k+1}^2 L C_n-\alpha_{k+1})f(y_k) \, +\\ + \, \E_{e_{k+1}}[V_{z_{k+1}}(u)\mid e_1, \ldots, e_{k}] - V_{z_{k}}(u) \leqslant \alpha_{k+1}f(u).
			\end{array}
		\end{equation}
	\end{lemma}
	
	Сложности возникают при попытке перенесения этого результата на случай $Q\ne \RR^n$. Ограничимся рассмотрением случая $p=q=2$, так как даже для него не удается обобщить рассуждения из \cite{Allen-Zhu-Orecchia-linear-coupling}. Введем обозначение
	\begin{equation*}
		\Prog_{s}(x) \overset{\text{def}}{=} -\underset{y\in Q}{\min}\left\{ \langle s, \, y-x \rangle + \frac{L}{2} \left \| y-x \right \|_2^2 \right\}.
	\end{equation*}
	Заметим, что вторая часть леммы 1 в евклидовом случае упрощается, а именно
	$$
	\E_e[\langle s,e\rangle^2] = \frac{||s||_2^2}{n},
	$$
	где $e\sim RS_2^n$ (см. лемму B.10 из \cite{pagerank}; формулировка данной леммы из \cite{pagerank} есть в Приложении~--- см. лемму~\ref{lemma_pagerank}). Отсюда следует, что $C_{n,q} = n^2$. Если положить
	$$
	\mbox{Grad}_e \left( x \right) \overset{\text{def}}{=} \underset{y\in Q}{\argmin}\left\{ \langle \langle\nabla f(x), \, e\rangle e, \, x-y \rangle - \frac{L}{2} \left \| y-x \right\|_2^2 \right\}
	$$
	и
	$$
	\mbox{Mirr}_e \left( {x,z,\alpha } \right)\overset{\text{def}}{=} \underset{y\in Q}{{\argmin}} \left\{ {\alpha \left\langle {n\left\langle {\nabla f\left( x \right), \, e} \right\rangle e, \, y-z} \right\rangle +V_z \left( y \right)} \right\},
	$$
	то чтобы обобщить приведенные рассуждения на условный случай, нужно оценить подходящим образом $\Prog_{n\langle\nabla f(x), \, e\rangle \, e}(x_{k+1})$ (точнее, его математическое ожидание по $e_{k+1}$), т.е., исходя из техники, используемой в \cite{Allen-Zhu-Orecchia-linear-coupling}, хотелось бы доказать оценку
	\begin{equation}\label{wrong_intuition}
		\E_{e_{k+1}}\left[\Prog_{n\langle\nabla f(x_{k+1}), \, e_{k+1}\rangle \, e_{k+1}}(x_{k+1})\right] \leqslant n^2 \left(f(x_{k+1})-\E \left [f \left (y_{k+1} \right ) \right ] \right ),
	\end{equation}
	чтобы получить оценку скорости сходимости, как и в случае безусловной минимизации.
	К сожалению, существует пример (будет приведен далее) выпуклой $L$-гладкой функции и замкнутого выпуклого множества, для которых \eqref{wrong_intuition} не выполняется.
	
     Сначала рассмотрим более детально $\Prog_{\xi}(x)$:
		\begin{equation*}
			\begin{array}{l}
				\Prog_{\xi}(x) = - \underset{y\in Q}{\min}\left\{\frac{L}{2}\Big|\Big| y-x \Big|\Big|_2^2 + \langle\xi, \, y-x\rangle \right\} =\\= - \underset{y\in Q}{\min}\left\{\frac{L}{2}\Big|\Big| y-x \Big|\Big|_2^2 + \langle\xi, \, y-x\rangle + \frac{1}{2L}\Big|\Big| \xi \Big|\Big|_2^2 \right\} + \frac{1}{2L}\Big|\Big| \xi \Big|\Big|_2^2=\\ 
				= - \underset{y\in Q}{\min}\left\{\Big|\Big|\frac{1}{\sqrt{2L}}\xi + \sqrt{\frac{L}{2}}\cdot y - \sqrt{\frac{L}{2}}\cdot x \Big|\Big|_2^2 \right\} + \frac{1}{2L}\Big|\Big| \xi \Big|\Big|_2^2=\\ = -\frac{L}{2}\underset{y\in Q}{\min}\left\{\Big|\Big| y - \left(x - \frac{1}{L}\xi\right) \, \Big|\Big|_2^2 \right\} + \frac{1}{2L}\Big|\Big| \xi \Big|\Big|_2^2,
			\end{array}
		\end{equation*}
		т.е. точка, в которой достигается этот минимум\footnote{См. \cite{Gasn_Universal_2018}.},
		$$
		\hat{y} = \pi_Q\left(x-\frac{1}{L}\xi\right).
		$$
		Тогда
		$$
		y_{k+1} = \pi_Q\left(x-\frac{1}{L}s_{k+1}\right), \quad s_{k+1} \overset{\text{def}}{=} \langle\nabla f(x_{k+1}), \, e_{k+1}\rangle \, e_{k+1} = \frac{1}{n}g(x_{k+1},e_{k+1}).
		$$
		Кроме того, обозначим через $\widetilde{y}_{k+1}$ точку множества $Q$, в которой достигается минимум в формуле для $\Prog_{n\langle\nabla f(x_{k+1}), \, e_{k+1}\rangle \, e_{k+1}}(x_{k+1})$. Тогда
		$$
		\widetilde{y}_{k+1} = \pi_Q\left(x-\frac{n}{L}s_{k+1}\right).
		$$
		Также для удобства рассмотрим следующие представления для $y_{k+1}$ и $\widetilde{y}_{k+1}$:
		\begin{equation}\label{example:representation_y_k+1}
		    y_{k+1} = x_{k+1} - \frac{1}{L}s_{k+1} + r_{k+1},
		\end{equation}
		$$
		\widetilde{y}_{k+1} = x_{k+1} - \frac{n}{L}s_{k+1} + \widetilde{r}_{k+1},
		$$
		где $r_{k+1}$ и $\widetilde{r}_{k+1}$ будем называть векторами невязок.

\begin{figure}[htb]
			\center{\includegraphics[width=0.65
				\linewidth]{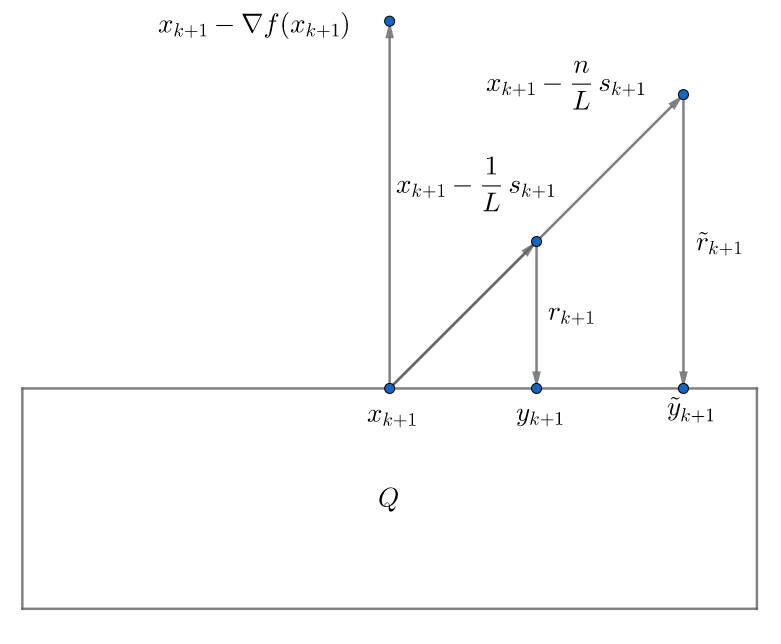}}
			\caption{Пример ситуации, когда ключевое неравенство не выполнено}
			\label{ris:bad_set}
		\end{figure}

		Рассмотрим функцию
		\begin{equation}\label{bad_function}
			f(y) = f(x_{k+1}) + \langle\nabla f(x_{k+1}), \, y-x_{k+1}\rangle + \frac{L}{2}\Big|\Big|y-x_{k+1}\Big|\Big|_2^2
		\end{equation}
		и множество, изображенное на рис.~\ref{ris:bad_set} (в качестве $\nabla f(x_{k+1})$ можно выбрать любой ненулевой вектор, а в качестве $Q$~--- прямоугольный параллелепипед с достаточно длинными сторонами, в центре одной из гиперграней которого размещена точка $x_{k+1}$).

		Подставим в \eqref{bad_function} значение $y=y_{k+1}$ и воспользуемся представлением $y_{k+1}$ из \eqref{example:representation_y_k+1}:
		$$
		-\langle \nabla f(x_{k+1}), \, -\frac{1}{L}s_{k+1} + r_{k+1} \rangle - \frac{L}{2}\Big|\Big|r_{k+1}-\frac{1}{L}s_{k+1}\Big|\Big|_2^2 = f(x_{k+1}) - f(y_{k+1}).
		$$
		Далее воспользуемся тем, что $s_{k+1} = \langle\nabla f(x_{k+1}), \, e_{k+1}\rangle \, e_{k+1}$:
		\begin{equation*}
		\begin{array}{rl}
		\frac{1}{L}\langle\nabla f(x_{k+1}), \, e_{k+1}\rangle^2 - \langle\nabla f(x_{k+1}), \, r_{k+1}\rangle - \frac{L}{2}\Big|\Big|r_{k+1}\Big|\Big|_2^2 + \langle r_{k+1}, \, s_{k+1}\rangle \, -\\
		- \, \frac{1}{2L}\langle\nabla f(x_{k+1}) , \, e_{k+1}\rangle^2 = f(x_{k+1}) - f(y_{k+1}),
		\end{array}
		\end{equation*}
		или в более компактной форме
		\begin{equation}\label{bad_config_key_inneq}
		\begin{array}{rl}
		\frac{1}{2L}\langle\nabla f(x_{k+1}), \, e_{k+1}\rangle^2 - \frac{L}{2}\Big|\Big|r_{k+1}\Big|\Big|_2^2 + \left \langle r_{k+1}, \,  s_{k+1} - \nabla f(x_{k+1}) \right \rangle=\\ = f(x_{k+1}) - f(y_{k+1}).
		\end{array}
		\end{equation}
		При таком выборе функции и множества получаем, что $n^2\cdot||r_{k+1}||_2^2 = ||\tilde{r}_{k+1}||_2^2$ для всех единичных $e$. Действительно, если
		$$
		\Prog_{\langle\nabla f(x_{k+1}), \, e_{k+1}\rangle \, e_{k+1}}(x_{k+1}) = \frac{1}{2L}\langle\nabla f(x_{k+1}) , \, e_{k+1}\rangle^2 - \frac{L}{2}\Big|\Big|r_{k+1}\Big|\Big|_2^2
		$$
		и
		$$
		\Prog_{n\langle\nabla f(x_{k+1}), \, e_{k+1}\rangle \, e_{k+1}}(x_{k+1}) = \frac{n^2}{2L}\langle\nabla f(x_{k+1}) , \, e_{k+1}\rangle^2 - \frac{L}{2}\Big|\Big|\tilde{r}_{k+1}\Big|\Big|_2^2,
		$$
		то
		$$
		\Prog_{n\langle\nabla f(x_{k+1}), \, e_{k+1}\rangle \, e_{k+1}}(x_{k+1}) = n^2 \, \Prog_{\langle\nabla f(x_{k+1}), \, e_{k+1}\rangle \, e_{k+1}}(x_{k+1}).
		$$
		Отсюда и из \eqref{bad_config_key_inneq} следует, что
		\begin{equation*}
		\begin{array}{rl}
		\frac{1}{n^2}\Prog_{n\langle\nabla f(x_{k+1}), \, e_{k+1}\rangle \,  e_{k+1}}(x_{k+1}) + \langle r_{k+1}, \, s_{k+1} - \nabla f(x_{k+1})\rangle=\\ = f(x_{k+1}) - f(y_{k+1}).
		\end{array}
		\end{equation*}
		Заметим, что вектор $s_{k+1}$ всегда короче (точнее, не длиннее) вектора $\nabla f(x_{k+1})$ и направлен ``вниз'' (т.е. в то же полупространство, образованное гранью $Q$, на которой лежит точка $x_{k+1}$), как и $\nabla f(x_{k+1})$. Значит, разность $s_{k+1} - \nabla f(x_{k+1})$ будет направлена в противоположную часть пространства. А вектор $r_{k+1}$ тоже направлен вниз. Следовательно, всегда выполняется $\langle r_{k+1}, \, s_{k+1} - \nabla f(x_{k+1})\rangle \leqslant 0$, причем с ненулевой вероятностью выполнено строгое неравенство. Это означает, что 
		$$
		\E_{e_{k+1}}\left[\langle r_{k+1}, \, s_{k+1} - \nabla f(x_{k+1})\rangle\right] < 0.
		$$
		Поэтому
		\begin{equation*}
			\begin{array}{l}
				\E_{e_{k+1}}\left[\Prog_{n\langle\nabla f(x_{k+1}), \, e_{k+1}\rangle \,  e_{k+1}}(x_{k+1})\right]=\\ = n^2(f(x_{k+1}) - \E_{e_{k+1}}[f(y_{k+1})]) - \E_{e_{k+1}}\left[\langle r_{k+1}, \, s_{k+1} - \nabla f(x_{k+1})\rangle\right]>\\ > n^2 \left ( f(x_{k+1}) - \E_{e_{k+1}} \left [f(y_{k+1}) \right ] \right ).
			\end{array}
		\end{equation*}

	Представленный контр-пример показывает трудности в перенесении предлагаемого в статье метода на задачи условной оптимизации. Теорема утверждает, что ускоренный неевклидов спуск (алгоритм ACDS) через $N$~итераций выдаст точку $y_N$, удовлетворяющую неравенству $\E[f(y_N)] - f(x^*) \leqslant \varepsilon,\, \varepsilon > 0$, если $N = O\left(\sqrt{\frac{\Theta L_2C_{n,q}}{\varepsilon}} \right)$. По определению $C_{n,q} = \sqrt{3}\min\{2q-1,32\ln n -8\}n^{\frac{2}{q}+1}$, а в случае $p= q =2$ можно взять $C_{n,q} = n^2$, что видно из леммы 1 и леммы B.10 из \cite{pagerank} (приведена в Приложении как лемма~\ref{lemma_pagerank}), поэтому $N = O\left(\sqrt{\frac{\Theta L_2n^2}{\varepsilon}} \right)$. Если же $p=1$ и $q=\infty$, то $C_{n,q} =n \sqrt{3} \, (32\ln n-8)$ и $N = O\left(\sqrt{\frac{\Theta L_2n\ln n}{\varepsilon}} \right)$.
	
	\section{Численные эксперименты}
	Для практического применения предложенный ускоренный неевклидов спуск по случайному направлению~ACDS был реализован на языке программирования~Python.
	Код метода и демонстрация
	вычислительных свойств метода
	с построением графиков сходимости доступны как~Jupyter Notebook и выложены
	в свободном доступе на~Github~\cite{git}.
	
	Рассмотрим следующую {\sl з\,а\,д\,а\,ч\,у}.
\textit{Пусть $A$~--- матрица размеров~$n \times n$
 со случайными независимыми элементами, равномерно распределенными на отрезке $[0, \, 1]$, а
матрица $B = \frac{A^\top A}{\lambda_{\max}(A^\top A)}$,
где $\lambda_{\max}(A^\top A)$~--- максимальное собственное значение матрицы $A^\top A$}.

Необходимо минимизировать функцию
\begin{equation}
\label{exp_func}
f = \frac{1}{2} \left \langle x - x_*, \, B(x - x_*) \right \rangle, 
\quad x \, \in \, \RR^n,
\end{equation}
где 
$x_* = (1, \, 0, \, 0, \, \ldots, \, 0)^\top$. 
Решение этой задачи известно и равно $x_*$, $f(x_*) = 0$. Начальная точка $x_0$ для всех экспериментов
выбиралась как $(0, \, 0, \, 0, \, \ldots, \, 1)^\top$.
Константа Липшица градиента целевой функции $L = 1$.

Следует отметить важность достаточно точного решения вспомогательной задачи 
	минимизации для нахождения $z_{k+1}$ на шаге~7 алгоритма~\ref{ACRCD} (зеркальный спуск).
В рассматриваемом случае эту задачу с помощью метода
множителей Лагранжа можно свести к 
задаче одномерной минимизации, подробности с формулами
см. в~\cite{git}. В реализации метода одномерная
минимизация выполняется с помощью обычной дихотомии
с точностью, на один порядок превышающей заданную.

			\begin{figure}[htb]
\center{\includegraphics[width=0.8\linewidth]{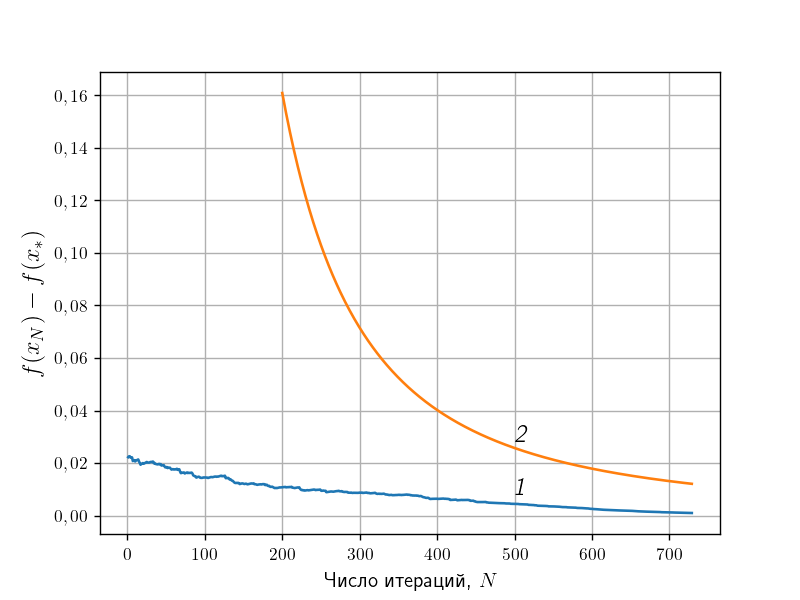}}
\caption{Сходимость ускоренного неевклидового спуска~ACDS для 
функции~\eqref{exp_func}, размерность $n = 10$. Показана практическая зависимость точности нахождения минимума $f(x_N)- f(x_*)$ от числа итераций~$N$ алгоритма (график~{\it 1}) и теоретическая оценка $O \left (\frac{4\Theta LC_{n,q}}{N^2} \right )$ (график~{\it 2})}
			\label{ris:exp10}
		\end{figure}

\begin{figure}[htb]
\center{\includegraphics[width=0.8\linewidth]{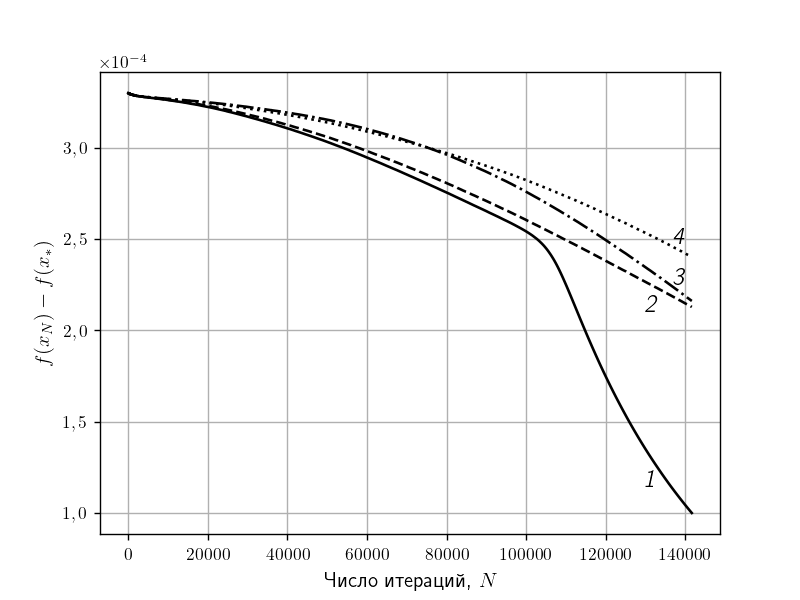}}
			\caption{Сходимость ускоренного неевклидового спуска~ACDS для 
функции~\eqref{exp_func}, размерность $n = 10^3$. Показана практическая зависимость точности нахождения минимума $f(x_N)- f(x_*)$ от числа итераций~$N$ алгоритма (сплошная линия~---
график~{\it 1}). Также для сравнения приведены
результаты работы метода
при других~$p$ (евклидова норма~--- график~{\it 2}; $\, p=1,8$~--- график~{\it 3}; $\, p=1,9$~--- график~{\it 4}) при тех же генерируемых векторах~$e$ и точке старта $x_0$}
			\label{ris:exp1000}
		\end{figure}

Для различных $n$ и заданной точности $\varepsilon$
были рассчитаны теоретически требуемые значения числа итераций
по теореме и проведена проверка сходимости на практике. Для данной задачи во всех случаях практическая
скорость сходимости по функции была выше.
Так, например, для $n = 10$ и $\varepsilon = 10^{-3}$
заданная точность была достигнута 
за 729~итераций (см. рис.~\ref{ris:exp10}), а теоретическая оценка
числа итераций дает  2537~итераций. 
Далее, для $n = 10^3$ и $\varepsilon = 10^{-4}$ 
теоретическая оценка
числа итераций дает не более чем~255972 итерации.
По факту алгоритм завершил работу
за~141643 итерации
(см. рис.~\ref{ris:exp1000}). 
Медленный спуск в начале
работы метода объясняется близостью начального
значения целевой функции к оптимальному именно
в данном примере ($f(x_0) = 0,00032983$). 

При проведении численных экспериментов
было обнаружено, что преимущество выбора прокс-структуры, связанной с $1$-нормой, возникает только
в пространствах от средней размерности (от $n=1000$). Будет ли иметь преимущество
предложенный метод, можно определить, сравнив
теоретические
оценки числа итераций предложенного метода
для разных $p$. На рис.~\ref{ris:exp1000} именно показан случай,
когда ускоренный спуск по направлению с неевклидовой прокс-структурой оказывается оптимальнее
спуска по направлению с евклидовой прокс-структурой.

В целом, численные эксперименты с ускоренным
спуском по случайному направлению подтверждают
теоретические результаты.

	\section{Заключение}
	В статье предложен ускоренный неевклидов спуск по направлению для решения задачи выпуклой безусловной оптимизации.
	В отличие от известных вариантов методов спуска по направлению (см., например \cite{nest}) в данной статье рассматривается ускоренный спуск по направлению с неевклидовой прокс-структурой. В случае когда $1$-норма решения близка к $2$-норме решения (это имеет место, например, если решение задачи разрежено~--- имеет много нулевых компонент), предлагаемый подход улучшает оценку на необходимое число итераций, полученную оптимальным методом из \cite{nest}, приблизительно в $\sqrt{n}$ раз, где $n$ -- размерность пространства, в котором происходит оптимизация.
	
	Данная статья открывает цикл работ, в которых планируется привести полные доказательства утверждений, полученных авторами в 2014--2016 гг. и приведенных (без доказательств) в \cite{disser}. В частности, далее планируется распространить приведенные в настоящей статье результаты на безградиентные методы, на задачи стохастической оптимизации и распространить все эти результаты на случай сильно выпуклой функции. 
	
	Открытой проблемой остается распространение полученных здесь результатов на случай задач оптимизации на множествах простой структуры. Напомним, что в статье существенным образом использовалось то, что оптимизация происходит на всем пространстве. Тем не менее в будущем планируется показать, что приведенные здесь результаты распространяются на задачи оптимизации на множествах простой структуры в случае, если градиент функционала в точке решения равен нулю (принцип Ферма).

    Авторы выражают благодарность Павлу Двуреченскому и Александру Тюрину за помощь в работе.

	\appendix{}

	\textsl{~~~Д~о~к~а~з~а~т~е~л~ь~с~т~в~о~\; л~е~м~м~ы~\;2.}\,
	Докажем сначала первую часть неравенства:
	\begin{equation}\label{lem2:main_estimation}
		\begin{array}{l}
			\alpha_{k+1}\langle n \langle \nabla f(x_{k+1}), \, e_{k+1}\rangle \, e_{k+1}, \, z_{k}-u\rangle =\\
			= \langle \alpha_{k+1}n \langle \nabla f(x_{k+1}), \, e_{k+1}\rangle \, e_{k+1}, \,  z_{k}-z_{k+1}\rangle \, +\\ 
			+ \, \langle \alpha_{k+1}n \langle \nabla f(x_{k+1}), \, e_{k+1}\rangle \, e_{k+1}, \,  z_{k+1}-u\rangle \overset{\circledOne}{\leqslant}\\
			\overset{\circledOne}{\leqslant} \langle \alpha_{k+1}n \langle \nabla f(x_{k+1}), \, e_{k+1}\rangle \, e_{k+1}, \,  z_{k}-z_{k+1}\rangle + \langle -\nabla V_{z_k}(z_{k+1}), \, z_{k+1}-u\rangle \overset{\circledTwo}{=}\\
			\overset{\circledTwo}{=} \langle \alpha_{k+1}n \langle \nabla f(x_{k+1}), \, e_{k+1}\rangle \, e_{k+1}, \, z_{k}-z_{k+1}\rangle + V_{z_k}(u) - V_{z_{k+1}}(u) - V_{z_k}(z_{k+1}) \overset{\circledThree}{\leqslant}\\
			\overset{\circledThree}{\leqslant} \left(\langle \alpha_{k+1}n \langle \nabla f(x_{k+1}), \, e_{k+1}\rangle \, e_{k+1}, \,  z_{k}-z_{k+1}\rangle - \frac{1}{2}||z_k-z_{k+1}||_p^2\right) +\\ 
			+ \, V_{z_k}(u) - V_{z_{k+1}}(u),
		\end{array}
	\end{equation}
	где $\circledOne$ выполнено в силу того, что $z_{k+1} = \argmin\limits_{z\in \RR^n}\left\{V_{z_k}(z) +  \alpha_{k+1}\langle n \langle \nabla f(x_{k+1}), \, e_{k+1}\rangle \, e_{k+1}, \, z\rangle\right\}$, откуда следует, что $\langle \nabla V_{z_k}(z_{k+1}) +  \alpha_{k+1} n \langle \nabla f(x_{k+1}), \, e_{k+1}\rangle \, e_{k+1}, \, u - z_{k+1}\rangle \geqslant 0$ для всех $u\in Q=\RR^n$, $\circledTwo$ выполнено в силу равенства треугольника для дивергенции Брэгмана~\footnote{Действительно, \begin{equation*}
			\begin{array}{rl}
				\forall x, \, y \in \RR^n\quad \langle -\nabla V_x(y), \, y-u\rangle = \langle\nabla d(x) - \nabla d(y), \, y-u\rangle = (d(u) - d(x) - \langle\nabla d(x), \, u-x\rangle ) \, -\\
				- \, (d(u) - d(y) - \langle\nabla d(y), \, u-y\rangle ) - (d(y) - d(x) - \langle\nabla d(x), \, y-x\rangle ) = V_x(u)-V_y(u)-V_x(y).
			\end{array}
		\end{equation*}}, $\, \circledThree$ выполнено, так как $V_x(y) \geqslant \frac{1}{2}||x-y||_p^2$ в силу сильной выпуклости прокс-функции $d(x)$.
		
		Теперь покажем, что
		\begin{equation*}
			\begin{array}{rl}
				\langle \alpha_{k+1}n \langle \nabla f(x_{k+1}), \, e_{k+1}\rangle \, e_{k+1}, \, z_{k}-z_{k+1}\rangle - \frac{1}{2}||z_k-z_{k+1}||_p^2\leqslant\\
				\leqslant \frac{\alpha_{k+1}^2n^2}{2} |\langle \nabla f(x_{k+1}), \, e_{k+1}\rangle|^2 \cdot ||e_{k+1}||_q^2.
			\end{array}
		\end{equation*}
		Действительно, в силу неравенства Гельдера
		\begin{equation*}
			\begin{array}{rl}
				\langle \alpha_{k+1}n \langle \nabla f(x_{k+1}), \, e_{k+1}\rangle \, e_{k+1}, \, z_{k}-z_{k+1}\rangle\leqslant\\
				\leqslant \alpha_{k+1}n ||\langle \nabla f(x_{k+1}), \, e_{k+1}\rangle \, e_{k+1}||_q\cdot||z_{k}-z_{k+1}||_p=\\
				= \alpha_{k+1}n|\langle \nabla f(x_{k+1}), \, e_{k+1}\rangle|\cdot||e_{k+1}||_q\cdot||z_{k}-z_{k+1}||_p,
			\end{array}
		\end{equation*}
		откуда
		\begin{equation}\label{technical_inequlity}
			\begin{array}{rl}
				\langle \alpha_{k+1} n \langle \nabla f(x_{k+1}), \, e_{k+1}\rangle \, e_{k+1}, \, z_{k}-z_{k+1}\rangle - \frac{1}{2}||z_k-z_{k+1}||_p^2\leqslant\\
				\leqslant \alpha_{k+1}n|\langle \nabla f(x_{k+1}), \, e_{k+1}\rangle|\cdot||e_{k+1}||_q\cdot||z_{k}-z_{k+1}||_p - \frac{1}{2}||z_k-z_{k+1}||_p^2.
			\end{array}
		\end{equation}
		Положим $t = ||z_k-z_{k+1}||_p$, $a = \frac{1}{2}$ и $b = \alpha_{k+1}n|\langle \nabla f(x_{k+1}), \, e_{k+1}\rangle|\cdot||e_{k+1}||_q$, тогда правая часть в \eqref{technical_inequlity} имеет вид
		$$
		bt-at^2.
		$$
		Если рассматривать полученное выражение как функцию от $t\in\RR$, то ее максимум при $t\in\RR$ равен (а значит, при $t\in\RR_+$ не превосходит) $\frac{b^2}{4a}$. Отсюда и из \eqref{technical_inequlity} следует неравенство
		\begin{equation}\label{lem2:golder}
			\begin{array}{rl}
				\langle \alpha_{k+1}n \langle \nabla f(x_{k+1}), \, e_{k+1}\rangle \, e_{k+1}, \, z_{k}-z_{k+1}\rangle - \frac{1}{2}||z_k-z_{k+1}||_p^2 \leqslant\\
				\leqslant \frac{\alpha_{k+1}^2n^2}{2} |\langle \nabla f(x_{k+1}), \, e_{k+1}\rangle|^2 \cdot ||e_{k+1}||_q^2.
			\end{array}
		\end{equation}
		
		Итак, учитывая \eqref{lem2:main_estimation} и \eqref{lem2:golder}, получаем, что
		\begin{equation*}
			\begin{array}{rl}
				\alpha_{k+1}\langle n \langle \nabla f(x_{k+1}), \, e_{k+1}\rangle \, e_{k+1}, \, z_{k}-u\rangle\leqslant\\
				\leqslant \frac{\alpha_{k+1}^2n^2}{2} |\langle \nabla f(x_{k+1}), \, e_{k+1}\rangle|^2 \cdot ||e_{k+1}||_q^2 + V_{z_k}(u) - V_{z_{k+1}}(u).
			\end{array}
		\end{equation*}
		Беря условное математическое ожидание $\E_{e_{k+1}}[\;\cdot\mid e_1, \, e_2, \, \ldots, \, e_k]$ от левой и правой частей последнего неравенства и пользуясь вторым неравенством из леммы~1, получаем, что
		\begin{equation*}
			\begin{array}{l}
				\alpha_{k+1}\langle \nabla f(x_{k+1}), \, z_{k}-u\rangle\leqslant\\ \leqslant
				\alpha_{k+1}^2\cdot\frac{C_{n,q}}{2n}||\nabla f(x_{k+1})||_2^2 + V_{z_k}(u) - \E_{e_{k+1}}[V_{z_{k+1}}(u)\mid e_1, \, e_2, \, \ldots, \, e_k],
			\end{array}
		\end{equation*}
		где $C_{n,q} {=} \sqrt{3}\min\{2q-1,32\ln n -8\}n^{\frac{2}{q}+1}$. Чтобы доказать вторую часть неравенства \eqref{grad_estimation}, покажем, что
		\begin{equation}\label{lem2:grad_estim}
			||\nabla f(x_{k+1})||_2^2 \leqslant 2nL\left(f(x_{k+1})-\E_{e_{k+1}} \left [f(y_{k+1})\mid e_1, \, e_2, \, \ldots, \, e_k \right ]\right).
		\end{equation}
		Во-первых, для всех $x, \, y\in\RR$
		\begin{equation*}
			\begin{array}{rl}
				f(y) - f(x) = \int\limits_{0}^{1}\langle \nabla f(x+\tau(y-x)), \, y-x\rangle d\tau=\\
				= \langle\nabla f(x), \, y-x\rangle + \int\limits_{0}^{1}\langle \nabla f(x+\tau(y-x))-\nabla f(x), \, y-x\rangle d\tau\leqslant\\
				\leqslant \langle\nabla f(x), \, y-x\rangle + \int\limits_{0}^{1}||\nabla f(x+\tau(y-x))-\nabla f(x)||_2\cdot||y-x||_2d\tau\leqslant\\
				\leqslant \langle\nabla f(x), \, y-x\rangle + \int\limits_{0}^{1}\tau L||y-x||_2\cdot||y-x||_2d\tau=\\
				= \langle\nabla f(x), \, y-x\rangle + \frac{L}{2}||y-x||_2^2,
			\end{array}
		\end{equation*}
		т.е.
		$$
		-\langle\nabla f(x), \, y-x\rangle - \frac{L}{2}||y-x||_2^2 \leqslant f(x) - f(y).
		$$
		Беря в последнем неравенстве $x = x_{k+1}$, $y=\text{Grad}_{e_{k+1}}(x_{k+1}) = x_{k+1} - \frac{1}{L} \langle \nabla f(x_{k+1}), \, e_{k+1}\rangle \, e_{k+1}$, получим, что
		\begin{equation*}
			\begin{array}{rl}
				f(x_{k+1}) - f(y_{k+1}) \geqslant \frac{1}{L} \langle \nabla f(x_{k+1}), \, e_{k+1}\rangle^2 - \frac{1}{2L}\langle \nabla f(x_{k+1}), \, e_{k+1}\rangle^2\cdot||e_{k+1}||_2^2=\\
				= \frac{1}{2L}\langle \nabla f(x_{k+1}), \, e_{k+1}\rangle^2.
			\end{array}
		\end{equation*}
		Возьмем от этого неравенства условное математическое ожидание $\E_{e_{k+1}}[\;\cdot\mid e_1, \, e_2, \, \ldots, \, e_k]$, используя лемму~B.10 из \cite{pagerank} (см. лемму~\ref{lemma_pagerank}), и получим неравенство~\eqref{lem2:grad_estim}. Лемма 2 доказана.
		
		\textsl{~~~Д~о~к~а~з~а~т~е~л~ь~с~т~в~о~\; л~е~м~м~ы~\;3.}\,
		Доказательство состоит в выписывании цепочки неравенств:
		\begin{equation*}
			\begin{array}{rl}
				\alpha_{k+1}(f(x_{k+1}) - f(u)) \leqslant\\
				\leqslant \alpha_{k+1}\langle\nabla f(x_{k+1}), \, x_{k+1}-u\rangle=\\
				= \alpha_{k+1}\langle\nabla f(x_{k+1}), \, x_{k+1}-z_k\rangle + \alpha_{k+1}\langle\nabla f(x_{k+1}), \, z_{k}-u\rangle\overset{\circledOne}{=}\\
				\overset{\circledOne}{=} \frac{(1-\tau_k)\alpha_{k+1}}{\tau_k}\langle\nabla f(x_{k+1}), \, y_{k}-x_{k+1}\rangle+\alpha_{k+1}\langle\nabla f(x_{k+1}), \, z_{k}-u\rangle\overset{\circledTwo}{\leqslant}\\
				\overset{\circledTwo}{\leqslant} \frac{(1-\tau_k)\alpha_{k+1}}{\tau_k} (f(y_k)-f(x_{k+1})) + \alpha_{k+1}\langle\nabla f(x_{k+1}), \, z_{k}-u\rangle\overset{\circledThree}{\leqslant}\\ \overset{\circledThree}{\leqslant} \frac{(1-\tau_k)\alpha_{k+1}}{\tau_k} (f(y_k)-f(x_{k+1})) \, + \\ 
				+ \, \alpha_{k+1}^2LC_{n,q}\cdot(f(x_{k+1}) - \E_{e_{k+1}}[f(y_{k+1})\mid e_1, \, e_2, \, \ldots, \, e_k]) \, +
				\\ 
				+ \, V_{z_k}(u) - \E_{e_{k+1}}[V_{z_{k+1}}(u)\mid e_1, \, e_2, \, \ldots, \, e_k]\overset{\circledFour}{=}\\
				\overset{\circledFour}{=} (\alpha_{k+1}^2LC_{n,q}-\alpha_{k+1})f(y_k) - \alpha_{k+1}^2LC_{n,q}\E_{e_{k+1}}[f(y_{k+1})\mid e_1, \, e_2, \, \ldots, \, e_k] \, +\\
				+ \, \alpha_{k+1}f(x_{k+1}) + V_{z_k}(u) - \E_{e_{k+1}}[V_{z_{k+1}}(u)\mid e_1, \, e_2, \, \ldots, \, e_k].
			\end{array}
		\end{equation*}
		Действительно, $\circledOne$ выполнено, так как $x_{k+1} \overset{\text{def}}{=} \tau_kz_k+(1-\tau_k)y_k\; \Leftrightarrow\; \tau_k(x_{k+1}-z_k) = (1-\tau_k)(y_k - x_{k+1})$, $\circledTwo$ следует из выпуклости $f(\cdot)$ и неравенства $1-\tau_k\geqslant0$, $\circledThree$ справедливо в силу леммы~2 и в $\circledFour$ используется равенство $\tau_k = \frac{1}{\alpha_{k+1}LC_{n,q}}$. Лемма 3 доказана.
		
		\textsl{~~~Д~о~к~а~з~а~т~е~л~ь~с~т~в~о~\; т~е~о~р~е~м~ы.}\,
	Заметим, что при $\alpha_{k+1} = \frac{k+2}{2LC_{n,q}}$ выполнено равенство
	$$
	\alpha_{k}^2LC_{n,q} = \alpha_{k+1}^2LC_{n,q} - \alpha_{k+1} + \frac{1}{4LC_{n,q}}.
	$$
	Возьмем для $k=0, \, 1, \, \ldots, \, N-1$ от каждого неравенства \eqref{coupling} леммы~3 математическое ожидание по $e_1, \, e_2, \, \ldots, \, e_N$, просуммируем полученные неравенства и получим
	$$
	\alpha_N^2LC_{n,q}\E[f(y_N)] + \sum\limits_{k=1}^{N-1}\frac{1}{4LC_{n,q}}\E[f(y_k)] + \E[V_{z_{N}}(u)] - V_{z_0}(u) \leqslant \sum\limits_{k=1}^{N}\alpha_kf(u),
	$$
	где $\E[\cdot] = \E_{e_1, \, e_2, \, \ldots, \, e_N}[\cdot]$.
	Положим $u=x^*$. Так как $\sum_{k=1}^{N}\alpha_k = \frac{N(N+3)}{4LC_{n,q}}$, $\E[f(y_k)] \geqslant f(x^*), V_{z_N}(u) \geqslant 0$ и $V_{z_0}(x^*) = V_{x_0}(x^*) \leqslant \Theta$, то выполняется неравенство
	$$
	\frac{(N+1)^2}{4LC_{n,q}}\E[f(y_N)] \leqslant \left(\frac{N(N+3)}{4LC_{n,q}} - \frac{N-1}{4LC_{n,q}}\right)f(x^*) + \Theta,
	$$
	откуда следует, что $\E[f(y_N)] \leqslant f(x^*) + \frac{4\Theta LC_{n,q}}{(N+1)^2}$. Теорема доказана.
	
	Приводим формулировку леммы B.10 из \cite{pagerank}. Отметим, что в доказательстве нигде не использовалось, что второй вектор в скалярном произведении (помимо $e$) есть градиент функции $f(x)$ (поэтому утверждение леммы~\ref{lemma_pagerank} остается верным для произвольного вектора $s\in\RR^n$ вместо $\nabla f(x)$).
	\begin{LemAp}\label{lemma_pagerank}
	    Пусть $e\sim RS_2^n(1)$ и вектор $s\in\RR^n$~--- некоторый вектор. Тогда
	    \begin{equation*}
	        \E_e[\langle s,\,e\rangle^2] = \frac{||s||_2^2}{n}.
	    \end{equation*}
	\end{LemAp}

\end{document}